\documentclass[letterpaper]{amsart}

\usepackage{style} 

\title{Classifying Primitive Solvable Permutation Groups of Rank 5 and 6}
\author{Anakin Dey}
\address{Anakin~Dey\\ Department of Mathematics \\ University of Illinois Urbana-Champaign \\
\href{mailto:anakind2@illinois.edu}
  {{\ttfamily\upshape anakind2@illinois.edu}}}

\author{Kolton O'Neal}
\address{Kolton~O'Neal\\ Department of Mathematics \\ University of Nebraska-Lincoln \\
\href{mailto:koneal6@huskers.unl.edu}
  {{\ttfamily\upshape koneal6@huskers.unl.edu}}}

\author{Duc Van Khanh Tran}
\address{Duc~Van~Khanh~Tran\\ Department of Mathematics \\ University of Texas at Austin \\
\href{mailto:duc.tranvk@utexas.edu}
  {{\ttfamily\upshape duc.tranvk@utexas.edu}}}

\author{Camron Upshur}
\address{Camron~Upshur\\ Department of Mathematics \\ Arizona State University \\
\href{mailto:cupshur@asu.edu}
  {{\ttfamily\upshape cupshur@asu.edu}}}

\author{Yong Yang}
\address{Yong~Yang\\ Department of Mathematics \\ Texas State University \\
  \href{mailto:yang@txstate.edu}
  {{\ttfamily\upshape yang@txstate.edu}}}


\begin{document}

\begin{abstract}
    Let $G$ be a finite solvable permutation group acting faithfully and primitively on a finite set $\Omega$. 
    Let $G_0$ be the stabilizer of a point $\alpha$ in $\Omega$.
    The rank of $G$ is defined as the number of orbits of $G_0$ in $\Omega$, including the trivial orbit $\{\alpha\}$. 
    In this paper, we completely classify the cases where $G$ has rank 5 and 6, continuing the previous works on classifying groups of rank 4 or lower.
\end{abstract}

\maketitle

\section{Introduction}\label{sec: intro}
Let $G$ be a primitive permutation group acting faithfully on a set $\Omega$, and let $G_0$ be the stabilizer of an element $\alpha \in \Omega$. 
The \emph{rank} of $G$ is defined as the number of orbits of $G_0 \acts \Omega$. 
Some work has been done to try to characterize these groups for arbitrary ranks~\cite{rank_of_prim_solv, finite_prim_solv}.
However, there has been significantly more progress for the classification of groups of low ranks.
Both solvable and insolvable groups of rank 2 have been classified~\cite{huppert_rank_2_solvable, hering_rank_2_insolvable}.
Various authors have studied primitive groups of rank 3~\cite{rank_3_higman1966, rank_3_higman1968, rank_3_kantor1982, rank_3_liebeck1987, rank_3_liebeck1986}.
Foulser completely classified primitive solvable groups of rank 3 and gave a partial classification of primitive solvable groups of rank 4 in~\cite{foulser}. 
Revitalization of this work has been made possible due to the advent of computer algebra systems such as GAP~\cite{GAP4}, making it more feasible to study the actions of finite groups~\cite{reg_orbit_4}.
Using GAP, Dolorfino~\etal\ completely classified primitive solvable groups of rank 4~\cite{prior_work_rank_4}. 
In this work, we classify primitive solvable permutation groups of rank 5 and 6. 
Note that while~\cite{prior_work_rank_4} relies on a coarse classification due to~\cite{foulser}, our paper is self-contained and does not use such a prior classification. 
In recent years, there has been further interest in obtaining stronger lower bounds for the number of conjugacy classes of a finite group~\cite{Nguyen2021, Keller2023, Gao2023}.
We believe that the results and methods in this paper may have applications in this direction.

If $G$ is a primitive solvable permutation group acting faithfully on $\Omega$, then $\abs{\Omega} = p^d$ for some prime $p$ and positive integer $d$~\cite{thatgermanpaper}. 
Moreover, $G$ has a minimal nontrivial normal subgroup $V$, which is an elementary abelian group such that $\abs{V} = p^d$, so $V$ behaves as a $d-$dimensional vector space over $\F_p$.
We can now decompose $G$ into $G = V \rtimes G_0$, where $G_0$ acts on $V$ as an irreducible subgroup of $\GL(V)$. 
Conversely, if an irreducible group $G_0$ acts on such a vector space $V$, then we can construct a primitive permutation group $G$ by taking the semidirect product $G = V \rtimes G_0$~\cite[Section 2]{pd_semi_direct_thingy}. 
Consequently, instead of classifying the rank of $G$ based on the number of orbits of $G_0$ in $\Omega$, we can equivalently consider the number of orbits of $G_0 \acts V$, as these values are equivalent~\cite[Definition 2.1]{foulser}.

Viewing $G_0$ as an irreducible subgroup of $\GL(V)$, we classify $G_0$ into three classes $\mathfrak{A}$, $\mathfrak{B}$, and $\mathfrak{L}$, following the standard analysis of solvable linear groups in~\cite{coxeter1965, foulser, huppert_rank_2_solvable, suprunenko1963}.
The class $\mathfrak{A}$ consists of primitive subgroups of $\Gamma(V)$, where $\Gamma(V)$ is the semilinear group of $V$ defined as
\[
    \Gamma(V) = \Set{x \to ax^\sigma | x \in \F_{q^n}, a \in \F_{q^n}^\times, \sigma \in \Gal\pqty{\F_{q^n} / \F_q}},
\]
and there are infinitely many groups in this class~\cite{reg_orbit_4}. 
$\mathfrak{B}$ contains the remaining primitive subgroups of $\GL(V)$, and $\mathfrak{L}$~contains the imprimitive subgroups of $\GL(V)$. 

Our main result is as follows.
\begin{thrm}\label{thrm: main_result}
    Suppose $G = V \rtimes G_0$ is a finite primitive solvable permutation group of rank at most 6, where $G_0$ acts on $V$ as an irreducible subgroup of $\GL(V)$. 
    At least one of the following is true:
    \begin{enumerate}
        \item $G_0 \in \mathfrak{A}\colon$ $G_0$ is a primitive subgroup of $\Gamma(V)$;
        \item $G_0 \in \mathfrak{B}\colon$ $G_0$ is one of the remaining primitive subgroups of $\GL(V)$ and appears in at least one row of \Cref{tab: final_results}; or
        \item $G_0 \in \mathfrak{L}\colon$ $G_0$ is an imprimitive subgroup of $\GL(V)$, and there exists an imprimitivity decomposition $V = \bigoplus\limits_{i = 1}^r V_i$.        
        \begin{enumerate}
            \item If $\rank(G) \leq 5$, then $r=2$, 3, or 4, and the action of $G_0$ on each $V_i \setminus \set{0}$ is transitive.
            \item If $\rank(G) = 6$, then $2 \leq r \leq 10$ and $r \neq 5$ where the action of $G_0$ on each $V_i \setminus \set{0}$ is transitive, or \mbox{$r = 2$ where the actions} of $G_0$ on $V_1 \setminus \set{0}$ and $V_2 \setminus \set{0}$ are doubly transitive~\cite[Proposition 2.5]{foulser}.
        \end{enumerate}
    \end{enumerate}
\end{thrm}

The structure of this paper is as follows. 
First, we describe the structure of primitive solvable permutation groups $G$ acting faithfully on a vector space $V$ in \Cref{sec: theory}.
Specifically, we do this for the case when $G_0 \in \mathfrak{B}$.
Then, in \Cref{sec: comp}, we leverage this structure in order to enumerate the groups in $\mathfrak{B}$ we are interested in using the computer algebra system GAP~\cite{GAP4}.
These groups are fully detailed in \Cref{sec: res}.
In particular, we provide examples illustrating the choice of wording in \Cref{thrm: main_result}~(2) of where we state $G_0 \in \mathfrak{B}$ appears in at least one row rather than exactly one row of \Cref{tab: final_results} (see \Cref{ex: mult_extra,ex: iso_diff_rank,,ex: iso_b1_b2}).
Finally, we pose possible directions for future work in \Cref{sec: future}.

\section{Preliminary Results}\label{sec: theory}

Let $G$ be a primitive solvable permutation group acting faithfully on a set $\Omega$.
As described before, this can be decomposed into a solvable subgroup $G_0 \leq \GL(V)$ acting irreducibly on a $d$-dimensional vector space $V$ over $\F_p$.
We now proceed to describe the structure when $G_0 \in \mathfrak{B}$.
Note that in this case, $G_0$ acts primitively on $V$.
It follows that $G_0$ also acts quasi-primitively, meaning that all non-trivial normal subgroups of $G_0$ act homogeneously on $V$.

\begin{thrm}[Theorem 2.2 of~\cite{reg_orbit_1}, Theorem 2.2 of~\cite{reg_orbit_2}, and Theorem 2.1 of~\cite{reg_orbit_3}]\label{thrm: structure_thrm}
    Suppose a finite solvable group $G_0$ acts faithfully, irreducibly, and quasi-primitively on a $d$-dimensional vector space $V$ over a finite field $\mathbb{F}$ of characteristic $p$. 
    Then, every normal abelian subgroup of $G_0$ is cyclic, and $G_0$ has normal subgroups $Z(E) \leq U \leq F \leq A \leq G_0$ and a characteristic subgroup $E \leq F$ such that the following statements hold: 
    \begin{enumerate}
        \item $F = EU$ is a central product where $Z(E) = E \cap U$;
        \item $F/U \cong E/Z(E)$ is a direct sum of completely reducible (or semi-simple) $G_0/F$-modules;
        \item There is a decomposition $E = E_1 \times \cdots \times E_s$ where each $E_i$ is an extraspecial $q_i$-group.
        We have that $\abs{E_i} = q_i^{2 m_i + 1}$ for some distinct primes $q_1, \ldots, q_s$ and some integer $m_i \geq 1$.
        Denote $e_i \defeq \sqrt{\abs{E_i / Z(E_i)}} = q_i^{m_i}$ and $e \defeq e_1 \cdots e_s$. 
        \mbox{We have that $e \mid d$ and $\gcd(p, e) = 1$};
        \item $A = C_{G_0}(U)$, and $A/F$ acts faithfully on $E/Z(E)$;
        \item $U$ is cyclic and acts fixed point freely on $W$, where $W$ is an irreducible submodule of $V_{U}$;
        \item $\abs{U}$ divides $p^k - 1$ for some $k \geq 1$, and $W$ can be identified with the span of $U$ which is isomorphic to $\mathrm{GF}\left(p^k\right)$; and
        \item $\abs{V} := n = \abs{W}^{eb}$ for some integer $b$.
    \end{enumerate}
\end{thrm}
Note that for $G_0 \in \mathfrak{B}$, we must have that $e > 1$, as it is the product of prime powers.

We now state a series of lemmas and propositions with the goal of restricting the various parameters stated in \Cref{thrm: structure_thrm} for the case that $\rank(G) \leq 6$.
The value $e$ as defined in \Cref{thrm: structure_thrm}~(3) is of great importance to our work.
Our first step is to set a bound on the value $e$ in the case of $\rank(G) \leq 6$. 
One immediate lower bound on the rank of $G$ is
\begin{equation}\label{eq: rank_bound}
    \rank(G) \geq \ceil{\frac{\abs{V} - 1}{\abs{G_0}}} + 1.
\end{equation}
This follows from the facts that at least one trivial orbit exists, that all orbits of $G_0 \acts V$ must partition $V$, and that the largest possible orbit size is $\abs{G_0}$.
We now build on \Cref{eq: rank_bound} further.
\begin{lem}[Lemma 2.3 of~\cite{reg_orbit_2} ]\label{lem: upper_bound_G}
    We have $\abs{G_0}$ divides $\dim(W) \cdot \abs{A/F} \cdot e^2 \cdot \pqty{\abs{W} - 1}$.
\end{lem}

\begin{prop}\label{prop: rank_lower_bound}
    Using the notation from \Cref{thrm: structure_thrm}, 
    \begin{equation}\label{eq: rank_bound_log}
        \rank(G) \geq \ceil{\frac{\abs{W}^{e} - 1}{\log_2\pqty{\abs{W}} \cdot \abs{A/F} \cdot e^2 \cdot (\abs{W} - 1)}} + 1.
    \end{equation}
\end{prop}

 \begin{pf}
    Combining \Cref{eq: rank_bound} with \Cref{thrm: structure_thrm}~(7) and \Cref{lem: upper_bound_G}, we get
    \begin{equation}\label{eq: rank_bound_for_W}
        \rank(G) \geq \ceil{\frac{\abs{V} - 1}{\abs{G_0}}} +1 \geq \ceil{\frac{\abs{W}^{eb} - 1}{\dim(W) \cdot \abs{A/F} \cdot e^2 \cdot (\abs{W} - 1)}} + 1.
    \end{equation}
    Then, note that $\abs{W} = p^k$ for some prime $p$ and integer $k$, so $\dim(W) \mid k$.
    Also, $k = \log_p\pqty{\abs{W}}$, which yields that \mbox{$\dim(W) \mid \log_p\pqty{\abs{W}}$}. 
    Thus, $\dim(W) \leq \log_p\pqty{\abs{W}}$.
    Since $p$ is a prime, we have that \mbox{$\dim(W) \leq \log_p\pqty{\abs{W}} \leq \log_2\pqty{\abs{W}}$}.
    Using this and \Cref{eq: rank_bound_for_W}, we get that
    
    \begin{align*}
        \rank(G) &\geq \ceil{\frac{\abs{W}^{eb} - 1}{\dim(W) \cdot \abs{A/F} \cdot e^2 \cdot (\abs{W} - 1)}} + 1 \\
                 &\geq \ceil{\frac{\abs{W}^{e} - 1}{\dim(W) \cdot \abs{A/F} \cdot e^2 \cdot (\abs{W} - 1)}} + 1 \\
                 &\geq \ceil{\frac{\abs{W}^{e} - 1}{\log_2\pqty{\abs{W}} \cdot \abs{A/F} \cdot e^2 \cdot (\abs{W} - 1)}} + 1.
    \end{align*}
\end{pf}

\begin{lem}\label{lem: AF_comp_red}
    $A / F$ acts completely reducibly on $E / Z(E)$.
\end{lem}
\begin{pf}
    By \Cref{thrm: structure_thrm} (2), $E/Z(E)$ is a direct sum of completely reducible $G_0/F$-modules.
    Also, by \Cref{thrm: structure_thrm}, $A \normal G_0$, so $E/Z(E)$ is a direct sum of completely reducible $A/F$-modules. 
    Therefore, $A/F$ acts completely reducibly on $E/Z(E)$.
\end{pf}

\begin{lem}\label{lem: alpha_lambda}
    Let $\alpha = \log_9\pqty{96 \cdot \sqrt[3]{3}}$ and $\lambda = 2 \cdot \sqrt[3]{3}$.
    Then, $\abs{A / F} \leq \frac{e^{2 \alpha}}{\lambda}$.
\end{lem}
\begin{pf}
    This is an immediate application of~\cite[Theorem 3.5]{manzwolf}.
    Note that $E/Z(E)$ is a direct sum of completely reducible $A/F$-modules by \Cref{lem: AF_comp_red} and that $A/F$ is solvable since $G$ is solvable. 
    Therefore, by~\cite[Theorem 3.5]{manzwolf}, 
    \[
        \abs{A/F} \leq \frac{\abs{E/Z(E)}^{\alpha}}{\lambda} = \frac{e^{2\alpha}}{\lambda}.
    \]
\end{pf}
\vspace{-10pt}
\begin{lem}\label{lem: lower_bound_w}
    Let $\rad(e)$ be the product of distinct prime factors of $e$. 
    Then, $\rad(e) \bigm| \pqty{\abs{W} - 1}$.
\end{lem}
\begin{pf}
    By \Cref{thrm: structure_thrm}~(5), $U$ acts fixed point freely on $W$. This implies that
    \[
        \abs{U} \bigm| \pqty{\abs{W} - 1}.
    \]
    By \Cref{thrm: structure_thrm}~(1) we have that $E\cap U = Z(E)$, so $Z(E) \leq U$. 
    Then, by \Cref{thrm: structure_thrm}~(3),
    \[
        \rad(e) \defeq q_1 \cdots q_s \bigm| \abs{Z(E)}.
    \]
    Therefore,
    \[
        \rad(e) \bigm| \abs{Z(E)} \bigm| \abs{U} \bigm| \pqty{\abs{W} - 1}.
    \]
\end{pf}
\begin{prop}\label{prop: e_values}
    Suppose $\rank(G) \leq 6$ and $G_0 \in \mathfrak{B}$.
    Then, the only possible values for $e$ are 2, 3, 4, 5, 6, 7, 8, 9, and 16.
\end{prop}
\begin{pf}
    Combining \Cref{eq: rank_bound_log} with \Cref{lem: alpha_lambda} yields
    \begin{equation}\label{eq: alpha_lambda_rank_bound}
        \rank(G) \geq \ceil{\frac{\abs{W}^{e} - 1}{\log_2\pqty{\abs{W}} \cdot \abs{A/F} \cdot e^2 \cdot (\abs{W} - 1)}} + 1 \geq \ceil{\frac{\lambda \cdot \pqty{\abs{W}^{e} - 1}}{\log_2\pqty{\abs{W}} \cdot e^{2 \alpha + 2} \cdot \pqty{\abs{W} - 1} }} + 1.
    \end{equation}
    By \Cref{lem: lower_bound_w} and the fact that the smallest possible value of $e$ is $2$, $\abs{W} \geq 3$. 
    One can see that the bound in \Cref{eq: alpha_lambda_rank_bound} is increasing with respect to $\abs{W}$ for $\abs{W} \geq 3$.
    By finding $e$ that satisfies the inequality
    \[
        \ceil{\frac{\lambda \cdot \pqty{3^{e} - 1}}{\log_2(3) \cdot e^{2 \alpha + 2} \cdot \pqty{3 - 1} }} + 1 > 6,
    \]
    we see that if $\rank(G) \leq 6$, then $e \leq 18$. We can now improve this bound on $e$ by considering individual cases. 
    
    \begin{table}[!ht]
        \begin{tabular}{|c|c|c|c|c|c|c|c|c|c|c|c|c|c|c|c|c|c|}
        \hline
        $e$   & 2 & 3 & 4 & 5 & 6 & 7 & 8 & 9 & 10 & 11 & 12 & 13 & 14 & 15 & 16 & 17 & 18 \\ 
        \hline
        $\rad(e)$ & 2 & 3 & 2 & 5 & 6 & 7 & 2 & 3 & 10 & 11 & 6  & 13 & 14 & 15 & 2  & 17 & 6  \\
        \hline
        $\abs{W} \geq$ & 3 & 4 & 3 & 11 & 7 & 8 & 3 & 4 & 11 & 23 & 7 & 27 & 29 & 16 & 3 & 103 & 7 \\
        \hline
        \end{tabular}
        \caption{Minimum values of $\abs{W}$ for each value of $e$.}\label{tab: w_lower_bound}
    \end{table}

    For each case, we find the minimum value of $\abs{W}$ by \Cref{lem: lower_bound_w} and the fact that $\abs{W}$ is a prime power by \Cref{thrm: structure_thrm}~(6).
    These minimum values of $\abs{W}$ are given in \Cref{tab: w_lower_bound}.
    If we take these values for $e$ and the corresponding lower bounds on $\abs{W}$ and apply them to \Cref{eq: alpha_lambda_rank_bound}, then we get that $\rank(G) > 6$ for $e \in \set{10, 11, 12,13,14,15,17,18}$.
    Therefore, if $\rank(G) \leq 6$, then $e \in \set{2, 3, 4, 5, 6, 7, 8, 9, 16}$ as claimed.
\end{pf}

Now that we have reduced the possible values of $e$ to a finite list, we aim to obtain upper bounds on $\abs{A / F}$ for each value of $e$.
By \Cref{eq: rank_bound_log}, larger values of $\abs{A / F}$ decrease the rank bound.
Thus, the largest values of $\abs{A / F}$ yield the worst case bound on $\rank(G)$.

\begin{lem}\label{lem: EZ_symp}
    $E / Z(E)$ is a symplectic vector space.
\end{lem}
\begin{pf}
    By~\cite[Corollary 1.10(iii)]{manzwolf},
    \[
        E / Z(E) = E_1 / Z(E) \times \cdots \times E_n / Z(E),
    \]
    where $E_i \leq C_{G_0}(E_j)$ for $i \neq j$. 
    In the proof of~\cite[Corollary 1.10]{manzwolf}, it is stated that $E_i/Z(E)$ has a non-degenerate symplectic form over $\F_{q_i}$. 
    Therefore, $E/Z(E)$ is a symplectic vector space.
\end{pf}

\begin{prop}\label{prop: AF_bounds}
    For each value of $e$ in \Cref{prop: e_values}, the following bounds for $\abs{A / F}$ hold.
    \begin{table}[!ht]
        \centering
        \begin{tabular}{|c|l|l|}
            \hline
            $e \defeq q^m$  & $\abs{A / F}$ divides one of the following values: & $\abs{A/F} \leq$\\
            \hline
            $2$             & $\Set{6}$     & $6$\\
            $3$             & $\Set{24}$    & $24$\\
            $4$             & $\Set{60, 6^2 \cdot 2}$ & $6^2\cdot 2$\\
            $5$             & $\Set{24}$ & $24$\\
            $6$             & $\Set{24 \cdot 6}$ & $24\cdot 6$\\
            $7$             & $\Set{48}$ & $48$\\
            $8$             & $\Set{42, 54, 6^3 \cdot 2, 6^4}$ & $6^4$\\
            $9$             & $\Set{40, 24^2 \cdot 2}$ & $24^2\cdot 2$\\
            $16$            & $\Set{136, 6^4 \cdot 24}$ & $6^4\cdot24$ \\
            \hline
        \end{tabular}
        \caption{Upper Bounds on $\abs{A/F}$}\label{tab: AF_upper_bound}
    \end{table}
\end{prop}
\begin{pf}
    For $e \in \set{2, 3, 4, 5, 7, 8, 9, 16}$, we let $e \defeq q^m$ for some prime $q$. 
    We have that $E / Z(E)$ is a symplectic vector space of dimension $2m$ over $\F_q$ by \Cref{lem: EZ_symp}, and $A / F$ acts completely reducibly on $E / Z(E)$ by \Cref{lem: AF_comp_red}. 
    Thus, we can apply~\cite[Lemma 2.17]{reg_orbit_1} with $G = A/F$ and $V = E/Z(E)$ to obtain the desired upper bounds for $\abs{A/F}$ when $e \in \set{2, 3, 4, 5, 7, 8, 9, 16}$. 
    Furthermore, the proof of~\cite[Lemma 2.17]{reg_orbit_1} describes many possible values that $\abs{A / F}$ divides for each value of $e$.
    
    For $e=6$, we let $q=6$ and $m=1$. In this case, we see from the proof of~\cite[Theorem 3.1]{reg_orbit_2} that we have $A / F \leq \SL(2, 3) \times \SL(2, 2)$, so $\abs{A / F} \bigm| 24 \cdot 6$.
\end{pf}

Using these values of $e$ and the corresponding bounds on $\abs{A / F}$, we can obtain upper bounds on $\abs{W}$.
\begin{prop}\label{prop: upper_bound_w}
    For each value of $e$ in \Cref{prop: e_values}, the following upper bounds on $\abs{W}$ hold when $\rank(G) \leq 6$.
    In particular, $e = 5$ and $e=7$ are ruled out when $\rank(G) \leq 6$.
    \begin{table}[ht]
        \centering
        \begin{tabular}{|c|c|c|c|c|c|c|c|}
            \hline
            $e$           & 2    & 3  & 4  & 6 & 8 & 9 & 16 \\
            \hline
            max $\abs{W}$ & 1511 & 79 & 31 & 7 & 7 & 4 & 3  \\
            \hline
        \end{tabular}
        \caption{Upper bounds for $\abs{W} = p^k$ for all possible values of $e$.}\label{tab: W_upper_bound}
    \end{table}
\end{prop}
\begin{pf}
    By \Cref{prop: rank_lower_bound}, $\rank(G) >6$ if
    \begin{equation}\label{eq: solve_w_for_6}
        \frac{\abs{W}^{e} - 1}{\log_2\pqty{\abs{W}} \cdot \abs{A/F} \cdot e^2 \cdot (\abs{W} - 1)} + 1 > 6.
    \end{equation}
    For each value of $e$, by using the largest corresponding bound of $\abs{A/F}$ from \Cref{prop: AF_bounds}, we can solve \Cref{eq: solve_w_for_6} for $\abs{W}$ to get an upper bound on $\abs{W}$ for which it is possible that $\rank(G) \leq 6$. 
    By \Cref{thrm: structure_thrm} (6) and \Cref{lem: lower_bound_w}, we can further improve the bound by restricting these values of $\abs{W}$ to be prime powers such that $\rad(e) \mid (\abs{W} - 1)$.
    This improvement completely eliminates the possibility of $e = 5$ and $e = 7$. 
    Therefore, we obtain upper bounds on $\abs{W}$ as listed in \Cref{tab: W_upper_bound}.
\end{pf} 

We now give an algorithm which optimizes the lower bound on $\rank(G)$ in the case that $e$ is a prime power. The first lower bound given in \Cref{prop: rank_lower_bound} was attained by using \Cref{eq: rank_bound}.
This bound is due to the fact that the largest possible size of an orbit of $G_0 \acts V$ is $\abs{G_0}$. 
In this algorithm, we extend this technique to account for the existence of smaller orbit sizes.

Fix some value $e = q^m$, where $q$ is prime, from \Cref{prop: e_values}.
Since $G_0$ has a trivial orbit in $V$, the other orbits contain a total of $\abs{V} - 1$ elements.
For each orbit $O$, $\abs{O}$ divides $\abs{G_0}$.
Let $B \defeq \dim(W) \cdot \abs{A/F}\cdot e^2 \cdot \pqty{\abs{W} - 1}$.
From \Cref{lem: upper_bound_G}, \mbox{$\abs{G_0}$ divides $B$}. 
Thus, $\abs{O}$ divides $B$, and the divisors of $B$ are possible sizes for the orbits of $G_0$. 

We now obtain specific parameters $q, m, p$, and $k$.
For each $e = q^m$, we get upper bounds on $\abs{W}$ and $\dim(W)$ over $\F_p$ from \Cref{prop: upper_bound_w}.
Let $p^k$ be some prime power less than or equal to this upper bound on $\abs{W}$, and the possible values for $\dim(W)$ are the different values of possible $k$.
We also get possible values of $\abs{A / F}$ to consider from the second column of \Cref{tab: AF_upper_bound} since we know that $\abs{A / F}$ divides one of those given values.
Note that we do not need to consider any proper divisors of the values in the second column of \Cref{tab: AF_upper_bound} since, by \Cref{eq: rank_bound_log}, smaller values of $\abs{A / F}$ can only make the rank bound greater.

The other values of consideration are $b$ and $d$, where $d$ is $\dim{(V)}$ over $\F_p$ such that $d = b \cdot k \cdot q^m$.
For each set of parameters $p$, $k$, $q$, and $m$, we check both of the cases where $b=1$ and where $b>1$. For a given set of parameters $p$, $k$, $q$, and $m$, we first let $b = 1$. 
Later, we describe how to consider $b > 1$.

\begin{alg}\label{alg: change}
Suppose we have parameters $p$, $k$, $q$, $m$, and $b=1$ as described as well as the corresponding values of $e$, $\abs{W}$, $\dim{(W)}$, and $\abs{A / F}$.
We know that for each orbit $O$, $\abs{O}$ divides $B \defeq \dim(W) \cdot \abs{A/F}\cdot e^2 \cdot \pqty{\abs{W} - 1}$.
Let $d_1,d_2,\ldots,d_t$ be the divisors of $B$, the possible sizes of these orbits.
In order to obtain a lower bound on the number of orbits of $G_0 \acts V$, we want to find the worst-case ``packing'' of elements of $V$ into orbits.
We know that there is at least one trivial orbit.
Thus, our goal is to pick $n_i$ orbits of size $d_i$ for each $d_i$ such that 
\[
    \sum_{i = 1}^t n_i d_i = \abs{V} - 1 = p^{b \cdot k \cdot q^m} - 1
\]
while minimizing
\[
    N = \sum_{i = 1}^t n_i.
\]
This is exactly the change-making problem, the problem of representing some chosen value using as few coins as possible from some fixed set of denominations~\cite{change_making}.
Let $N$ be the optimal solution to the change-making problem with the chosen value $\abs{V} - 1$ and the denominations being the set of divisors $\set{d_1, \ldots, d_t}$.
Finding $N$ is quite simple and can be done via a standard dynamic programming algorithm.
Then, the value $N + 1$ sets a lower bound on the number of orbits of $G_0 \acts V$, \ie\ a lower bound on the rank of $G$.
\end{alg}

\Cref{alg: change} described above helps restrict possible sets of parameters $e = q^m$ and $\abs{W} = p^k$ when $e$ is a prime power. When $e=6$, we let $q=6$ and $m=1$, and we have $p=7$ and $k=1$ from \Cref{prop: upper_bound_w}. 
We now describe how we consider $b > 1$.
For each of these sets of parameters, we use \Cref{eq: rank_bound_for_W} to bound the value of $b$. 
For a given set of parameters $e = q^m$, $\abs{W} = p^k$, and $b \geq 1$, if \Cref{eq: rank_bound_for_W} gives a bound less than or equal to 6 for any of the values of $\abs{A / F}$ from \Cref{prop: AF_bounds}, we keep that set of parameters. 
Otherwise, we eliminate that $b$ value and higher $b$ values for the corresponding set of $q$, $m$, $p$, and $k$, as the rank bound given by \Cref{eq: rank_bound_for_W} increases as $b$ increases.

Hence, \Cref{alg: change} and \Cref{eq: rank_bound_for_W} yield possible sets of parameters $q$, $m$, $p$, $k$, and $b$ for families of groups $G_0 \in \mathfrak{B}$ which have 6 or fewer orbits when acting on $V$.
We further classify $G_0 \in \mathfrak{B}$ into two cases based on if $b = 1$ or $b > 1$, denoted as $\mathfrak{B}_1$ and $\mathfrak{B}_>$, respectively.
It is due to the algorithmic process in \Cref{alg: main} that we split $\mathfrak{B}$ into these two cases rather than based on the action $E \acts V$ as in~\cite{prior_work_rank_4} or the action of a minimal normal abelian subgroup of $G_0$ acting on $V$ as in~\cite{foulser}.

\begin{thrm}
    Let $G = V \rtimes G_0$, where $G_0 \in \mathfrak{B}$, and $V$ is a vector space over a field of characteristic $p$.
    Then, by \Cref{thrm: structure_thrm}, $G_0$ has an extraspecial subgroup $E$ of order $q^{2m + 1}$ for some prime $q$, and $q^m = e$.
    Let $d = \dim(V)$ and $k$ the integer guaranteed by \Cref{thrm: structure_thrm}~(6) such that $d = b\cdot k \cdot q^m$. 
    If $G_0$ is not described by the sets of parameters in \Cref{tab: b_1_params,tab: b_2_params}, then $\rank(G) > 6$.
    In addition, the case where $e = 6$ is also included. 
\end{thrm}

\begin{table}[H]
    \small
    \centering
    \begin{tabular}[t]{llllllc}
        \toprule
        No. & $q$ & $m$ & $p$ & $k$ & $d$ & Rank $\geq$ \\ 
        \midrule
        1 & 2 & 1 & 3 & 1 & 2 & 2 \\ 
        2 & 2 & 1 & 3 & 2 & 4 & 3 \\ 
        3 & 2 & 1 & 3 & 3 & 6 & 3 \\ 
        4 & 2 & 1 & 3 & 4 & 8 & 4 \\ 
        5 & 2 & 1 & 3 & 5 & 10 & 4 \\ 
        6 & 2 & 1 & 5 & 1 & 2 & 2 \\ 
        7 & 2 & 1 & 5 & 2 & 4 & 3 \\ 
        8 & 2 & 1 & 5 & 3 & 6 & 4 \\ 
        9 & 2 & 1 & 7 & 1 & 2 & 2 \\ 
        10 & 2 & 1 & 7 & 2 & 4 & 3 \\ 
        11 & 2 & 1 & 11 & 1 & 2 & 2 \\ 
        12 & 2 & 1 & 11 & 2 & 4 & 5 \\ 
        13 & 2 & 1 & 13 & 1 & 2 & 3 \\ 
        14 & 2 & 1 & 13 & 2 & 4 & 6 \\ 
        15 & 2 & 1 & 17 & 1 & 2 & 3 \\ 
        16 & 2 & 1 & 19 & 1 & 2 & 3 \\ 
        17 & 2 & 1 & 23 & 1 & 2 & 2 \\ 
        18 & 2 & 1 & 29 & 1 & 2 & 3 \\ 
        \bottomrule
    \end{tabular}
    \hspace{10pt}
    \begin{tabular}[t]{llllllc}
        \toprule
        No. & $q$ & $m$ & $p$ & $k$ & $d$ & Rank $\geq$ \\ 
        \midrule
        19 & 2 & 1 & 31 & 1 & 2 & 3 \\ 
        20 & 2 & 1 & 37 & 1 & 2 & 4 \\ 
        21 & 2 & 1 & 41 & 1 & 2 & 4 \\ 
        22 & 2 & 1 & 43 & 1 & 2 & 4 \\ 
        23 & 2 & 1 & 47 & 1 & 2 & 3 \\ 
        24 & 2 & 1 & 53 & 1 & 2 & 4 \\ 
        25 & 2 & 1 & 59 & 1 & 2 & 4 \\ 
        26 & 2 & 1 & 61 & 1 & 2 & 5 \\ 
        27 & 2 & 1 & 67 & 1 & 2 & 5 \\ 
        28 & 2 & 1 & 71 & 1 & 2 & 4 \\ 
        29 & 2 & 1 & 73 & 1 & 2 & 5 \\ 
        30 & 2 & 1 & 79 & 1 & 2 & 5 \\ 
        31 & 2 & 1 & 83 & 1 & 2 & 5 \\ 
        32 & 2 & 1 & 89 & 1 & 2 & 6 \\ 
        33 & 2 & 1 & 97 & 1 & 2 & 6 \\ 
        34 & 2 & 1 & 101 & 1 & 2 & 6 \\ 
        35 & 2 & 1 & 103 & 1 & 2 & 6 \\ 
        36 & 2 & 1 & 107 & 1 & 2 & 6 \\ 
        \bottomrule
    \end{tabular}
    \hspace{10pt}
    \begin{tabular}[t]{llllllc}
        \toprule
        No. & $q$ & $m$ & $p$ & $k$ & $d$ & Rank $\geq$ \\ 
        \midrule
        37 & 2 & 2 & 3 & 1 & 4 & 2 \\ 
        38 & 2 & 2 & 3 & 2 & 8 & 4 \\ 
        39 & 2 & 2 & 5 & 1 & 4 & 3 \\ 
        40 & 2 & 2 & 7 & 1 & 4 & 3 \\ 
        41 & 2 & 2 & 11 & 1 & 4 & 4 \\ 
        42 & 2 & 2 & 13 & 1 & 4 & 5 \\ 
        43 & 2 & 3 & 3 & 1 & 8 & 4 \\ 
        44 & 2 & 3 & 5 & 1 & 8 & 5 \\ 
        45 & 3 & 1 & 2 & 2 & 6 & 3 \\ 
        46 & 3 & 1 & 2 & 4 & 12 & 4 \\ 
        47 & 3 & 1 & 5 & 2 & 6 & 4 \\ 
        48 & 3 & 1 & 7 & 1 & 3 & 3 \\ 
        49 & 3 & 1 & 13 & 1 & 3 & 4 \\ 
        50 & 3 & 1 & 19 & 1 & 3 & 5 \\ 
        51 & 3 & 2 & 2 & 2 & 18 & 6 \\ 
        52 & 6 & 1 & 7 & 1 & 6 & 2 \\[22pt]
        \bottomrule
    \end{tabular}
    \vspace{-7pt}
    \caption{Parameters for $G_0 \in \mathfrak{B}_1$ which possibly have 6 or fewer orbits on $V$.}\label{tab: b_1_params}
    \vspace{5pt}
    \begin{tabular}{llllllc}
        \toprule
        No. & $q$ & $m$ & $p$ & $k$ & $d$ & Rank $\geq$ \\ 
        \midrule
        53 & 2 & 1 & 3 & 1 & 4 & 4 \\ 
        54 & 2 & 2 & 3 & 1 & 8 & 6 \\ 
        \bottomrule
    \end{tabular}
    \vspace{-7pt}
    \caption{Parameters for $G_0 \in \mathfrak{B}_>$ which possibly have 6 or fewer orbits on $V$.}\label{tab: b_2_params}
    \vspace{-7pt}
\end{table}

\begin{pf}
    \Cref{tab: b_1_params} and \Cref{tab: b_2_params} result from computing the minimum rank of all possible sets of parameters allowed by \Cref{prop: upper_bound_w} using \Cref{alg: change} and keeping only those whose minimum ranks are below 6.
\end{pf}

\section{Computation}\label{sec: comp}

Using the parameters in \Cref{tab: b_1_params,tab: b_2_params}, we can enumerate possible groups $G_0 \acts V$ with $6$ or fewer orbits using GAP~\cite{GAP4}.
The procedure used mirrors that of~\cite{prior_work_rank_4, reg_orbit_4}, and our computational process is split into two cases.

\begin{alg}\label{alg: main}
    When $G_0 \in \mathfrak{B}_1$, we take the following 4 steps:
    \begin{enumerate}
        \item[1:] Construct the extraspecial group $E$ as guaranteed by \Cref{thrm: structure_thrm}~(3) as well as its normalizer $N_E$ as subgroups of $\GL\pqty{e, p^k}$.
        \item[2:] Embed $N_E$ into $\GL\pqty{k \cdot e, p}$ using the embedding $\GL\pqty{e, p^k} \into \GL\pqty{k \cdot e, p}$.
        \item[3:] Construct the normalizer $N$ of $N_E$ in $\GL\pqty{k \cdot e, p}$.
        \item[4:] Since $G_0$ is a subgroup of $N$, enumerate all subgroups $G_0$ of $N$ up to conjugacy in $\GL\pqty{k \cdot e, p}$ and check for primitive solvable groups of 6 or fewer orbits which contain $E$ as an isomorphic subgroup of $G_0$.
    \end{enumerate}
    
    For $G_0 \in \mathfrak{B}_>$, we introduce an intermediate step after step 1:
    \begin{enumerate}
        \item[1.5:] Embed $N_E$ from $\GL\pqty{e, p^k}$ into $\GL\pqty{b \cdot e, p^k}$ by taking the tensor product of $N_E$ with $\GL\pqty{b, p^k}$.
        Call this new product $N_E$ and replace $k \cdot e$ with $b \cdot k \cdot e$ as needed for the rest of the algorithm.
    \end{enumerate}
\end{alg}

Notice that $e = 6$ is the only value of $e$ which is not a prime power.
Since $6$ is the product of two distinct primes, we can take advantage of a decomposition lemma.

\begin{lem}[Lemma 2.1(6) of~\cite{e6_decomp}]\label{lem: e6_decomp}
    Let $G_0 \acts V$ be a primitive solvable subgroup of $\GL(6, 7)$.
    Then, $G_0$ is conjugate in $\GL(6, 7)$ to the Kronecker product $G_2 \mathbin{\dot\times} G_3$, where $G_2$ is a primitive solvable subgroup of $\GL(2, 7)$, and $G_3$ is a primitive solvable subgroup of $\GL(3, 7)$.
\end{lem}

Our problem for $e=6$ is then reduced to enumerating primitive solvable subgroups of $\GL(2, 7)$ and $\GL(3, 7)$ and taking their Kronecker product by \Cref{lem: e6_decomp}.
This is a less computationally intensive process than \Cref{alg: main} in the case of $e = 6$.
However, this cannot be extended to the cases of $e$ being a prime power since it relies on the fact that $6$ is the product of distinct primes.

\section{Results}\label{sec: res}

\Cref{tab: final_results} lists all lines from \Cref{tab: b_1_params,tab: b_2_params} that correspond to families of groups $G = V \rtimes G_0$ such that $\rank(G) \leq 6$.
Recall that there are exactly two types of extraspecial groups of order $q^{2m + 1}$ for a prime $q$ and integer $m \geq 1$~\cite{robinson_group_theory}.
One of them is of exponent $q$, and the other is of exponent $q^2$. 
We denote these as extraspecial types, or et, $+$ and $-$, respectively.
For each set of parameters $q, m, p, k, d, b$, et we describe the number of groups $G_0$ of rank $2$ through $6$ with these parameters.
However, we note that there are some caveats to this list.

\begin{ex}\label{ex: mult_extra}
For a given $q$ and $m$, we may have that a group $G_0$ contains both extraspecial groups of type~$+$ and of type~$-$.
One such example of this is line~3 of \Cref{tab: b_1_params}.
When taking these parameters and the extraspecial group of type~$+$, one of the generated groups $G_0$ is $\Z_{13} \times \text{QD}_{16}$.
This group is also generated when taking the same parameters and the extraspecial group of type~$-$.
We count the group in both sets of parameters in this case.
\end{ex}

\begin{ex}\label{ex: iso_diff_rank}
We can have that a group is of multiple ranks.
What this means is that two isomorphic copies of a group $G_0$ may be in separate conjugacy classes of the larger group $N$ constructed in \Cref{alg: main} and act on $V$ differently, resulting in different ranks.
One case where this occurs is line~48 of \Cref{tab: b_1_params}.
When considering the extraspecial group of type~$+$, one of the rank $4$ groups and one of the rank $5$ groups are isomorphic to $\pqty{\pqty{\Z_3 \times \Z_3} \rtimes \Z_3} \rtimes Q_8$.
We also count the group in both ranks in this case.
\end{ex}

\begin{ex}\label{ex: iso_b1_b2}
We may have a group in $\mathfrak{B}_1$ that also appears in $\mathfrak{B}_>$.
Furthermore, in each of these cases, the rank may even be the same.
One such pair of parameters that demonstrates this is line~2 of \Cref{tab: b_1_params} with an extraspecial group of type~$-$ and line~53 of \Cref{tab: b_2_params} with an extraspecial group of type~$-$.
Both of these sets of parameters generate a group isomorphic to $\pqty{\pqty{\Z_2 \times \Z_2 \times \Z_2} \rtimes \pqty{\Z_2 \times \Z_2}} \rtimes \Z_3$.
Furthermore, both of these groups are of rank $4$.
Such examples are also considered in~\cite{prior_work_rank_4}.
Again, we count such groups in both sets of parameters.
\end{ex}

\begin{table}
    \small
    \centering
    \begin{tabular}{llllllllccccc}
        \toprule
        No. & $q$ & $m$ & $p$ & $k$ & $d$ & $b$ & et & \# Rank 2 & \# Rank 3 & \# Rank 4 & \# Rank 5 & \# Rank 6 \\
        \midrule
        1 & 2 & 1 & 3 & 1 & 2 & 1 & $-$ & 4 & 0 & 0 & 0 & 0\\ 
        1 & 2 & 1 & 3 & 1 & 2 & 1 & $+$ & 1 & 0 & 0 & 0 & 0\\ 
        2 & 2 & 1 & 3 & 2 & 4 & 1 & $-$ & 0 & 8 & 3 & 0 & 0\\ 
        3 & 2 & 1 & 3 & 3 & 6 & 1 & $-$ & 0 & 3 & 3 & 1 & 1\\ 
        3 & 2 & 1 & 3 & 3 & 6 & 1 & $+$ & 0 & 1 & 0 & 1 & 0\\ 
        4 & 2 & 1 & 3 & 4 & 8 & 1 & $-$ & 0 & 0 & 0 & 3 & 5\\ 
        5 & 2 & 1 & 3 & 5 & 10 & 1 & $-$ & 0 & 0 & 1 & 0 & 1\\ 
        6 & 2 & 1 & 5 & 1 & 2 & 1 & $-$ & 3 & 0 & 0 & 0 & 0\\ 
        7 & 2 & 1 & 5 & 2 & 4 & 1 & $-$ & 0 & 0 & 4 & 5 & 6\\ 
        8 & 2 & 1 & 5 & 3 & 6 & 1 & $-$ & 0 & 0 & 0 & 1 & 0\\ 
        9 & 2 & 1 & 7 & 1 & 2 & 1 & $-$ & 3 & 3 & 1 & 1 & 0\\ 
        9 & 2 & 1 & 7 & 1 & 2 & 1 & $+$ & 1 & 1 & 1 & 0 & 0\\ 
        10 & 2 & 1 & 7 & 2 & 4 & 1 & $-$ & 0 & 0 & 0 & 3 & 3\\ 
        11 & 2 & 1 & 11 & 1 & 2 & 1 & $-$ & 2 & 1 & 1 & 0 & 2\\ 
        11 & 2 & 1 & 11 & 1 & 2 & 1 & $+$ & 0 & 1 & 0 & 0 & 0\\ 
        13 & 2 & 1 & 13 & 1 & 2 & 1 & $-$ & 0 & 1 & 1 & 1 & 2\\ 
        15 & 2 & 1 & 17 & 1 & 2 & 1 & $-$ & 0 & 3 & 0 & 2 & 0\\ 
        16 & 2 & 1 & 19 & 1 & 2 & 1 & $-$ & 0 & 1 & 2 & 0 & 2\\ 
        16 & 2 & 1 & 19 & 1 & 2 & 1 & $+$ & 0 & 0 & 1 & 0 & 0\\ 
        17 & 2 & 1 & 23 & 1 & 2 & 1 & $-$ & 1 & 1 & 1 & 0 & 0\\ 
        17 & 2 & 1 & 23 & 1 & 2 & 1 & $+$ & 0 & 1 & 0 & 1 & 0\\ 
        18 & 2 & 1 & 29 & 1 & 2 & 1 & $-$ & 0 & 1 & 1 & 0 & 1\\ 
        19 & 2 & 1 & 31 & 1 & 2 & 1 & $-$ & 0 & 1 & 0 & 3 & 0\\ 
        19 & 2 & 1 & 31 & 1 & 2 & 1 & $+$ & 0 & 0 & 1 & 0 & 1\\ 
        20 & 2 & 1 & 37 & 1 & 2 & 1 & $-$ & 0 & 0 & 1 & 0 & 1\\ 
        21 & 2 & 1 & 41 & 1 & 2 & 1 & $-$ & 0 & 0 & 1 & 1 & 2\\ 
        22 & 2 & 1 & 43 & 1 & 2 & 1 & $-$ & 0 & 0 & 1 & 0 & 1\\ 
        23 & 2 & 1 & 47 & 1 & 2 & 1 & $-$ & 0 & 1 & 0 & 1 & 0\\ 
        23 & 2 & 1 & 47 & 1 & 2 & 1 & $+$ & 0 & 0 & 0 & 1 & 0\\ 
        24 & 2 & 1 & 53 & 1 & 2 & 1 & $-$ & 0 & 0 & 1 & 0 & 1\\ 
        25 & 2 & 1 & 59 & 1 & 2 & 1 & $-$ & 0 & 0 & 1 & 0 & 1\\ 
        26 & 2 & 1 & 61 & 1 & 2 & 1 & $-$ & 0 & 0 & 0 & 1 & 0\\ 
        27 & 2 & 1 & 67 & 1 & 2 & 1 & $-$ & 0 & 0 & 0 & 1 & 0\\ 
        28 & 2 & 1 & 71 & 1 & 2 & 1 & $-$ & 0 & 0 & 1 & 0 & 0\\ 
        28 & 2 & 1 & 71 & 1 & 2 & 1 & $+$ & 0 & 0 & 0 & 0 & 1\\ 
        29 & 2 & 1 & 73 & 1 & 2 & 1 & $-$ & 0 & 0 & 0 & 0 & 1\\ 
        30 & 2 & 1 & 79 & 1 & 2 & 1 & $-$ & 0 & 0 & 0 & 1 & 0\\ 
        31 & 2 & 1 & 83 & 1 & 2 & 1 & $-$ & 0 & 0 & 0 & 1 & 0\\ 
        32 & 2 & 1 & 89 & 1 & 2 & 1 & $-$ & 0 & 0 & 0 & 0 & 1\\ 
        34 & 2 & 1 & 101 & 1 & 2 & 1 & $-$ & 0 & 0 & 0 & 0 & 1\\ 
        35 & 2 & 1 & 103 & 1 & 2 & 1 & $-$ & 0 & 0 & 0 & 0 & 1\\ 
        36 & 2 & 1 & 107 & 1 & 2 & 1 & $-$ & 0 & 0 & 0 & 0 & 1\\ 
        37 & 2 & 2 & 3 & 1 & 4 & 1 & $-$ & 3 & 4 & 0 & 0 & 0\\ 
        37 & 2 & 2 & 3 & 1 & 4 & 1 & $+$ & 0 & 7 & 6 & 0 & 0\\ 
        39 & 2 & 2 & 5 & 1 & 4 & 1 & $-$ & 0 & 0 & 0 & 1 & 4\\ 
        39 & 2 & 2 & 5 & 1 & 4 & 1 & $+$ & 0 & 0 & 5 & 0 & 8\\ 
        40 & 2 & 2 & 7 & 1 & 4 & 1 & $-$ & 0 & 1 & 0 & 1 & 1\\ 
        40 & 2 & 2 & 7 & 1 & 4 & 1 & $+$ & 0 & 0 & 0 & 0 & 3\\ 
        43 & 2 & 3 & 3 & 1 & 8 & 1 & $-$ & 0 & 0 & 0 & 5 & 5\\ 
        43 & 2 & 3 & 3 & 1 & 8 & 1 & $+$ & 0 & 0 & 0 & 2 & 3\\ 
        45 & 3 & 1 & 2 & 2 & 6 & 1 & $+$ & 0 & 5 & 2 & 0 & 0\\ 
        46 & 3 & 1 & 2 & 4 & 12 & 1 & $+$ & 0 & 0 & 0 & 1 & 3\\ 
        48 & 3 & 1 & 7 & 1 & 3 & 1 & $+$ & 0 & 0 & 3 & 2 & 2\\ 
        49 & 3 & 1 & 13 & 1 & 3 & 1 & $+$ & 0 & 0 & 0 & 0 & 1\\ 
        53 & 2 & 1 & 3 & 1 & 4 & 2 & $-$ & 0 & 10 & 5 & 0 & 0\\ 
        53 & 2 & 1 & 3 & 1 & 4 & 2 & $+$ & 0 & 8 & 3 & 0 & 0\\ 
        54 & 2 & 2 & 3 & 1 & 8 & 2 & $-$ & 0 & 0 & 0 & 2 & 0\\ 
        \bottomrule
    \end{tabular}
    \vspace{5pt}
    \caption{All parameters of $G_0 \in \mathfrak{B}$ which describe families of groups of rank $6$ or lower.}\label{tab: final_results}
\end{table}

It is because of these examples that \Cref{thrm: main_result} states that a given $G_0 \in \mathfrak{B}$ appears in at least one row of \Cref{tab: final_results} rather than exactly one row.

We provide all families of $G_0$ for the parameters described in \Cref{tab: final_results} as separate GAP files.
These files, along with the details for the implementation of the algorithms in \Cref{sec: comp} in GAP, can be found on \href{https://github.com/Spamakin/Solvable-Primitive-Permutation-Groups-of-Rank-5-and-6}{GitHub}\footnote{\href{https://github.com/Spamakin/Solvable-Primitive-Permutation-Groups-of-Rank-5-and-6}{https://github.com/Spamakin/Solvable-Primitive-Permutation-Groups-of-Rank-5-and-6}}.
This repository also contains notes on naive optimizations used to make \Cref{alg: main} faster as well as use less memory.
Furthermore, it also contains an implementation of \Cref{alg: change} and files explicitly describing Examples \ref{ex: mult_extra}, \ref{ex: iso_diff_rank}, and \ref{ex: iso_b1_b2}.

We also note some minor corrections to the classification of primitive solvable permutation groups of rank 4 due to~\cite{prior_work_rank_4}.
This prior work found only 3 distinct groups for line 37 in Table 4 in the case of groups of rank at most 4.
However, we find 4 distinct groups.
Furthermore, the case found in~\cite{prior_work_rank_4} with $b > 1$ and $k = 1$ as stated has parameters which define a family of groups of rank significantly larger than $4$.
Instead, that case must have $b = 1$ and $k = 5$, which then corresponds to line~5 in \Cref{tab: final_results}.

\section{Future Work}\label{sec: future}
The methodology of this paper can be extended in a straightforward manner to classify groups of higher rank. 
The analysis in \Cref{prop: e_values} can be extended by considering the relevant inequalities with a value higher than 6. 
One could apply the analysis of \Cref{sec: theory}, replacing 6 as desired with some larger rank, and then list possible parameters for higher ranks in a similar manner to \Cref{tab: b_1_params} and \Cref{tab: b_2_params}. 
Starting with \Cref{prop: e_values}, we may consider setting \Cref{eq: alpha_lambda_rank_bound} to be greater than or equal to some higher value than 6 and proceed with similar analysis from there obtaining bounds on the parameters $p, k, q$ and $m$.
The same computational method described in \Cref{sec: comp} can then be applied. 

The main source of difficulty of extending our method to higher ranks such as 7 and 8 comes from the computational intensity of constructing groups with larger parameters for higher ranks.
Steps 1 and 4 of \Cref{alg: main} require us to solve the subgroup isomorphism problem, and step 4 requires us to enumerate subgroups.
Both of these problems are some of the hardest problems in computational group theory~\cite{Gary1978, Pyber1996}.
As one considers larger and larger ranks, the sizes of our vector spaces $V$, and thus their general linear groups, grow at a considerable rate.
As a concrete example, Line 51 in \Cref{tab: b_1_params} took 8 days of computation to verify that no groups of rank $6$ or below exist for that particular set of parameters.

Further analysis on the structure of primitive solvable groups is required to overcome this difficulty.
While this issue of extending to higher ranks could be solved with further computational power, this avenue becomes prohibitively expensive.
Rather, future work should focus on theoretical improvements to obtaining bounds on the parameters $q$, $m$, $p$, $k$, and $d$ as well as optimizations to \Cref{alg: main}.

\section*{Acknowledgements}
This research was conducted under NSF-REU grant DMS-1757233 and DMS-2150205 by Anakin Dey, Kolton O'Neal, Duc Van Khanh Tran, and Camron Upshur during the Summer of 2023 under the supervision of Prof.\ Yong Yang.
The authors gratefully acknowledge the financial support of NSF and thank Texas State University for providing a great working environment and support.
The authors would also like to thank Professors Thomas Keller and Derek Holt for their invaluable help.
Finally, we would like to thank the reviewer for their remarks which have improved this manuscript.

\section*{Disclosure Statement}
The authors declare that there is no potential conflict of interest that could influence this paper.

\section*{Data Availability Statement} Data sharing is not applicable to this article, as no data sets were generated or analysed during the current study.

\printbibliography

\end{document}